\documentstyle[amssymb,12pt]{article}

\begin{document}

\title{Functionals on Closed 2-Surfaces}

\vspace{5cm}

\maketitle

\begin{center}
Metin G{\"{u}}rses\\
 {\small Department of Mathematics, Faculty of Sciences, Bilkent University}\\
  {\small 06800 Ankara, Turkey, gurses@fen.bilkent.edu.tr}

\end{center}

\begin{abstract}
We show that the 2-torus in ${\mathbb R}^3$ is a critical point of a sequence of functionals ${\cal F}_{n}$ ($n=1,2,3, \cdots$) defined over compact 2-surfaces in ${\mathbb R}^3$. When the Lagrange function ${\cal E}$ is a polynomial of degree $n$ of the mean curvature $H$ of the surface, the radii ($a,r$) of the 2-torus are related as  $\frac{a^2}{r^2}=\frac{n^2-n}{n^2-n-1},~~ n \ge 2$. If the Lagrange function  depends on both mean and Gaussian curvatures, the 2- torus remains to be a critical point of ${\cal F}_{n}$ without any  constraints on the radii of the torus.
\end{abstract}

\section{Introduction}

In the history of differential geometry
there are some special subclasses of $2$-surfaces, such as surfaces of
constant Gaussian curvature, surfaces of constant
mean curvature, minimal surfaces and the Willmore surfaces. These surfaces arise in many
different branches of sciences; in particular, in various parts of
theoretical physics (string theory, general theory of relativity),
cell-biology and differential geometry \cite{part}-\cite{qu}. All these special surfaces constitute
critical points of certain functionals. Although the Euler-Lagrange equations of these functionals are very
complicated and difficult, there are some techniques to solve these equations. Using the
deformation of the Lax equations of the integrable equations it is possible to construct
2-surfaces in ${\mathbb R}^3$  solving more general shape equations \cite{sym1}-\cite{gurses}.

The main objective in our work is to investigate 2-surfaces  derivable from a variational principle,
such as the minimal surfaces, Willmore surfaces
\cite{will1}, \cite{will2} and surfaces solving the shape equation
\cite{tu1}-\cite{Yang-can}. In all these cases the Lagrange function is a quadratic function of the mean
curvature of the surface. It is natural to ask whether there are
surfaces solving the Euler-Lagrange equations corresponding to
more general Lagrange functions depending on the mean and Gaussian
curvatures of the surface \cite{tu1}-\cite{tu3}. It is the  purpose of this work to find an answer to such a
question.

The quadratic Helfrich functional \cite{hel1973} for a theoretical model of a closed cell-membrane is
\begin{equation}\label{hel}
{\cal F}=\frac{1}{2}\,\int_{S}\,[k_{c}\,(2\,H+c_{0})^2+2\,w]\, dA+p\, \int_{V}\, dV
\end{equation}
where $k_{c}$ is the elasticity constant, $H$ and $c_{0}$ are the mean and the spontaneous curvatures,  $w$ is the surface tension and $p$ is the pressure difference between in and out of the surface. First variation of the above functional gives the shape equation

\begin{equation}\label{shape}
p-2wH+2\,k_{c}\, \nabla^2\, H+k_{c}\,(2H+c_{0})(2H^2-c_{0}H-2K)=0
\end{equation}
Sphere with an arbitrary radius is an exact solution  and shape equation (\ref{shape}) gives a relation between the
model parameters and the radius of the sphere. Stability of this solution has been studied in \cite{hel1}. A special 2-torus, called the Clifford torus, is also an exact solution of (\ref{shape}) \cite{tu1}. We might consider that Helfrich's functional (\ref{hel})  as  an approximation of the higher order functionals ${\cal F}_{n}=\int_{S}\, {\cal E}_{n}\, dA+p \int_{V}\, dV$, ~~($n \ge 2$) where

\begin{eqnarray}
{\cal E}_{n}&=&\sum_{k=0}^{n}\, a_{n+1-k}\, H^{k},~~n=1,2, \cdots \nonumber\\
&=&a_{1}\,H^{n}+a_{2}\,H^{n-1}+\cdots+a_{n}\,H+a_{n+1},
\end{eqnarray}
Here $a_{i}$'s are constants describing the parameters of the cell-membrane model. Hence it is worthwhile to study such functionals and search for possible critical points. This is another motivation of this work on functionals on closed surfaces.

 In section 2 we introduce the the first and second variations of a functional ${\cal F}=\int_{S}\, {\cal E} dA+p \int_{V}\, dV$ where
 the Lagrange function is an arbitrary function of the mean and Gauss curvatures of the surface $S$ \cite{tu1}.  In section 3 we consider the Lagrange function ${\cal E}$ depends only on the  the mean curvature. Letting ${\cal E}_{n}$ be a polynomial of degree $n$, then we show that
 2-torus $T_{n}$ is a critical point of the corresponding functional ${\cal F}_{n}=\int_{S}\, {\cal E}_{n} dA+p \int_{V}\, dV$ with
$\frac{a^2}{r^2}=\frac{n^2-n}{n^2-n-1},~~ n=2,3,4, \cdots $. When $n=2$ the corresponding torus $T_{2}$ is the Clifford torus. In section 4 we consider ${\cal E}_{n}$ depend on both mean and Gaussian curvatures. We show that, inclusion of suitable Gaussian curvature terms into ${\cal E}_{n}$, 2-torus with arbitrary radii $(a,r)$ becomes a critical point. In section 5 we prove two theorems: First one is to prove
 the constraint equation for the radii of the 2-torus and the second one is on the more general functionals without the constraint equation. In Appendix A and B we give more examples of functionals discussed in sections 3 and 4.

\section{Surfaces from a variational principle}

Let $S$ be a $2$-surface in ${\mathbb R}^3$ with Gaussian  ($K$) and mean ($H$) curvatures.  A functional
${\cal F}$ is defined by

\begin{equation}\label{funct}
{\cal F}=\int_{S}\, {\cal E}(H,K) dA+p \int_{V} dV
\end{equation}

\noindent where ${\cal E}$ is the Lagrange function depending on $H$ and $K$. Functional ${\cal F}$ is
also called curvature energy or shape energy. Here $p$ is
a constant which play the role of Lagrange multiplier  and $V$ is the volume enclosed within the surface $S$.
We obtain the Euler-Lagrange equations corresponding to the above functional from its first variation.
Let ${\cal E}$ be a twice differentiable function of $H$ and $K$. Then the first variation of ${\cal F}$ is given by
\begin{equation}\label{first}
\delta\,{\cal F}=\int_{S}\, E({\cal E}) \, \Omega\, dA.
\end{equation}
where $\Omega$ is an arbitrary smooth function on $S$.
Then the Euler-Lagrange
equation $E({\cal E})=0$ for ${\cal F}$ reduces to \cite{tu1}-\cite{tu3}

\begin{equation}\label{el1}
(\nabla^2 +4H^2-2K) {\partial {\cal E} \over \partial
H}+2({\nabla} \cdot \bar{\nabla}+2 K H) {\partial {\cal E} \over
\partial K}-4H {\cal E} +2p=0.
\end{equation}
Here and in what follows, $\nabla^2={\frac{1}{\sqrt{g}}}{\frac{\partial}{\partial x^i}}
 \left({\sqrt{g}}{g^{ij}}{\frac{\partial}{\partial
 x^j}}\right)$and $\nabla{\cdot}{\bar{\nabla}}={\frac{1}{\sqrt{g}}}{\frac{\partial}{\partial x^i}}
 \left({\sqrt{g}}K{h^{ij}}{\frac{\partial}{\partial
 x^j}}\right)$, ${g}={\det{(g_{ij})}}$, $g^{ij}$ and $h^{ij}$ are
 inverse components of the first and second fundamental forms;
 $x^{i}=(u,v)$ and we assume Einstein's summation convention on
 repeated indices over their ranges.

For the second variation of the functional we assume that ${\cal E}$ depends only on $H$. In this case the
expression is much simpler \cite{tu1}

\begin{eqnarray}\label{second}
\delta^2\, {\cal F}&=&\int_{S}\,\left({\cal E}_{1}\, \Omega^2\, +{\cal E}_{2}\, \Omega\, \nabla^2\, \Omega-2 \frac{\partial {\cal E}}{\partial H}\, \Omega \nabla \cdot \tilde{\nabla}\, \Omega +\frac{1}{4}\, \frac{\partial^2 {\cal E}}{\partial H^2}\, (\nabla^2 \Omega)^2 \right.\nonumber \\
&&\left.+\frac{\partial {\cal E}}{\partial H}\, [\nabla (H \Omega) \cdot \nabla \Omega-\nabla \Omega \cdot \tilde{\nabla} \Omega] \right) \,dA
\end{eqnarray}
where
\begin{eqnarray}
{\cal E}_{1}&=& (2H^2-K)^2\, \frac{\partial^2 {\cal E}}{\partial H^2}-2HK \frac{\partial {\cal E}}{\partial H}+2K {\cal E}-2 H p, \nonumber\\
{\cal E}_{2}&=& (2H^2-K)\, \frac{\partial^2 {\cal E}}{\partial H^2}+2H \frac{\partial {\cal E}}{\partial H}-{\cal E},
\end{eqnarray}
To have minimal energy solutions of (\ref{el1}) it is expected that the second variation $\delta^2 {\cal F} >0$.

\noindent
We have the following classical examples:
\begin{itemize}
    \item [{\bf i)}]Minimal surfaces: $ {\cal E}=1, ~~p=0$.
    \item [{\bf ii)}]Constant mean curvature surfaces: ${\cal E}=1$.
    \item [{\bf iii)}]Linear Weingarten surfaces: ${\cal E}=aH+b$, where $a$
and $b$ are some constants.
    \item [{\bf iv)}]Willmore surfaces: $ {\cal E}=H^2$ \cite{will1}, \cite{will2}.
    \item [{\bf v)}]Surfaces solving the shape equation (\ref{shape}) of lipid bilayer cell membranes: ${\cal E}=\frac{1}{2}\, k_{c}\,(2H+c_{0})^2+w$,
where $k_{c}$, $c_{0}$ and $w$ are  constants \cite{hel1}-\cite{tu3}.
\end{itemize}

2-sphere  in ${\mathbb R}^3$ has constant mean and Gaussian curvatures. Hence it is
 a critical point of the most  general functional (\ref{funct}). Eq. (\ref{el1}) gives a relation
 between $p$, radius of the sphere and other parameters in model.

 Another compact surface in ${\mathbb R}^3$ is the 2-torus, $T$. It has been shown \cite{Yang-can} that a special kind of torus, known as the Clifford torus, solves the shape equation (example (v) above). In this work we shall show that, $T$ solves not only a quadratic shape equation but it solves all Euler-Lagrange equations with Lagrange function ${\cal E}$  is a polynomial function of the mean and Gaussian curvatures. 2-Torus in ${\mathbb R}^3$ is defined as $X: U \to {\mathbb R}^3 $ where \cite{do}
\begin{eqnarray}
&&X(u,v) = ((a + r \cos u) \cos v, (a + r \cos u) \sin v, r \sin u),\\
&&0 < u < 2\pi, 0 < v > 2\pi
\end{eqnarray}
The first and second fundamental forms are

\begin{eqnarray}
ds_{1}^2=g_{ij}\,dx^{i}\,dx^{j}= r^2\,du^2 + (a + r \cos u)^2\,dv^2, \label{tor1}\\
ds_{2}^2=h_{ij}\,dx^{i}\,dx^{j}= r^2\,du^2 + (a + r \cos u) \cos u dv^2 \label{tor2}
\end{eqnarray}
The Gaussian ($K$) and mean ($H$) curvatures are
\begin{equation}\label{tor3}
K =\frac{\cos u}{r(a + r \cos u)},~~ H =\frac{1}{2}\, \left(\frac{1}{r}+\frac{\cos u}{(a + r \cos u)}\right)
\end{equation}
where $a$ and $r$  ($a > r$) are the radii of the torus. $K$ and $H$
satisfy a linear equation $r^2 K-2 r H+1=0$. Hence torus is a linear Weingarten surface.

\section{Functionals with Mean Curvature}

In this section we shall consider the Lagrange function ${\cal E}$ depending only on the mean
curvature $H$  of the surface. Furthermore se shall assume that ${\cal E}$ is a polynomial function of
$H$. Let the degree of the polynomial be $n$, then we write
\begin{eqnarray}
{\cal E}_{n}&=&\sum_{k=0}^{n}\, a_{n+1-k}\, H^{k},~~n=1,2, \cdots \nonumber\\
&=&a_{1}\,H^{n}+a_{2}\,H^{n-1}+\cdots+a_{n}\,H+a_{n+1},
\end{eqnarray}
where $a_{k}$, $(k=1,2, \cdots)$ are constants to be determined. Assuming that the torus is a critical point of the functional ${\cal F}$ we shall determine the coefficients $a_{i}$ of the polynomial expansion of ${\cal E}$ in terms of the torus radii $a$ and $r$. We shall give three examples here in this section and three more in Appendix A.

\vspace{0.5cm}
\noindent
1. {\bf First Order Functional}:  ${\cal E}_{1}=a_{1}\,H+a_{2}$

\vspace{0.5cm}
\noindent
{\bf Solution:}

\begin{equation}
p=-\frac{a_{1}}{r^2},~~a_{2}=-\frac{a_{1}}{r}
\end{equation}
There is no restriction on the radii $a$ and $r$. Torus is a critical point of the functional (\ref{el1}) for all values of $r$ and $a$.

\vspace{0.5cm}
\noindent
2. {\bf Second Order Functional}:  Most General case: ${\cal E}_{2}=a_{1}\,H^2+a_{2}\,H+a_{3}+a_{4}\,K$.
Here the term containing $a_{4}$ is a topologically invariant and gives no contribution to shape equation (\ref{el1}).

\vspace{0.5cm}
\noindent
{\bf Solution}:
\begin{equation}
p=-\frac{a_{2}}{r^2},~~a_{3}=-\frac{a_{2}}{r},~~a^2=2 r^2
\end{equation}
We find that (for $p=0$)
\begin{equation}\label{min-en2}
{\cal F}_{2}=\int_{T}\, {\cal E}_{2} dA=2 \pi^2\, a_{1}
\end{equation}
2-torus with $a^2=2 r^2$ minimizes the functional $\int_{S}\, {\cal E}_{2} dS$ and the minimum energy is given
in (\ref{min-en2}). Willmore conjecture \cite{will1},\cite{will-conj} states that ($a_{1} \ne 0$)
\begin{equation}
\frac{1}{a_{1}}\,\int_{S}\, {\cal E}_{2} dA \ge 2 \pi^2\,
\end{equation}
for all compact surfaces $S$ with genus $g>0$. The proof of this conjecture has been given very recently \cite{min-max}.
In terms of the Helfrich's functional (\ref{hel}) we have $a_{1}=k_{c}$, $a_{2}=2 k_{c} c_{0}$ and $a_{}=\frac{1}{2}\, k_{c} c_{0}^2+w$.
The parameters $p,w, c_{0}$ and $p$ must satisfy
\begin{equation}
p=\frac{2 k_{c} c_{0}}{r^2},~~w=p\,r\,(1+\frac{1}{4} r c_{0})
\end{equation}

Since $K$ is a topological invariant and  total curvature for torus is zero, there will be no contribution of single $K$ terms in ${\cal E}$. Because of this reason will not include a linear $K$  term in the ${\cal E}_{n}$ expressions.

\vspace{0.5cm}
\noindent
3. {\bf Third Order Functional}:  ${\cal E}_{3}=a_{1}\,H^3+a_{2}\,H^2+a_{3}\,H+a_{4}$.

\vspace{0.5cm}
\noindent
{\bf Solution}:
\begin{eqnarray}
p&=&-\frac{3a_{1}-a_{3} r^2}{r^4},~~a_{4}=\frac{2 a_{1}- a_{3} r^2}{ r^3},~~a_{2}=\frac{15 a_{1}}{2r}, \\
a^2&=&(6/5) r^2
\end{eqnarray}
We also find that (for $p=0$)
\begin{equation}\label{min-en3}
{\cal F}_{3}=\int_{T}\, {\cal E}_{3} dA=9 \sqrt{5}\, \pi^2\, (a_{1}/r)
\end{equation}
2-torus with $a^2=(6/5) r^2$ minimizes the functional $\int_{S}\, {\cal E}_{3} dS$ and the minimum energy is given
in (\ref{min-en3}). We expect that ($a_{1} \ne 0$)
\begin{equation}
(a_{1}/r)^{-1}\,\int_{S}\, {\cal E}_{3} dA \ge 9 \sqrt{5}\, \pi^2\,
\end{equation}
for all compact surfaces $S$ with genus $g>0$. This is the Willmore conjecture for $n=3$.

\vspace{0.5cm}
\noindent

It is possible to continue on this study and find infinitely many functionals with Lagrange functions ${\cal E}_{n}$  which are polynomials of degree $n$. We added three more examples of ${\cal E}_{n}$ with $n=4,5,6$ in Appendix A.
We observe that whatever the polynomial degree of the ${\cal E}_{n}$ is, there exists a constraint on the radii $a$ and $r$ of the torus which depend on the degree of the polynomials

\begin{equation}\label{const}
\frac{a^2}{r^2}=\frac{n^2-n}{n^2-n-1},~~ n=2,3,4, \cdots
\end{equation}
Hence for all these surfaces $ 1 < \frac{a^2}{r^2} \le 2$. The proof of this constraint equation is given in section 5.

These constraints can be utilized to select the correct functional for the toroidal fluid membranes. These functionals are
used to minimize the energy of the lipid membranes. The ratio $a/r$ of toroidal configuration can be measured experimentally. Comparing  the measured value of this
ratio with   those listed above (\ref{const}) we can find the degree of the polynomial function ${\cal E}_{n}$ from (\ref{const}), hence identifying the functional for the corresponding closed membrane. As an example, for vesicle membranes such a measurement had been done by Mutz and Bensimon \cite{M-B}. They measured the value of this ratio approximately as $\frac{a}{r}=1.43$, or $\frac{a^2}{r^2}=2.04$. Hence for vesicle membranes the correct functional should be the quadratic one which was first introduced by W. Helfrich several years ago \cite{hel1973}. For other closed fluid membranes the functionals might be different. In his study of vesicles with toroidal topology Seifert \cite{seifert1}-\cite{wilson} claimed that circular toroidal configurations exist for any $0 < v < 1$ where $v$ is the reduced volume parameter defined by
\begin{equation}
v=\frac{V}{\frac{4 \pi}{3} (\frac{A}{4 \pi})^{3/2}}
\end{equation}
Here $V$ is the volume and $A$ is the surface are of the toroidal membranes. In terms of the  2-torus radii we have
\begin{equation}
\frac{a^2}{r^2}=\frac{1}{1.94 v^4}
\end{equation}
 Seifert claims also that the membranes with $v \ge 0.84$ are not stable. This corresponds to $\frac{a^2}{r^2} \le 1.188$. This stability analysis seems to save the quadratic and cubic functionals we have given above. One must be careful about discarding higher order functionals, because
 the stability analysis mentioned by Seifert is done by using  the quadratic functional. To investigate the stability of toroidal configurations with $\frac{a^2}{r^2}=\frac{n^2-n}{n^2-n-1},~~ n=2,3,4, \cdots$ one must use the  corresponding functional ${\cal F}_{n}$, $(n=3,4,5,6, \cdots)$.

\section{Functionals with Mean and Gaussian Curvatures }

\noindent
In this section we shall consider the functionals which are composed of both curvatures $K$ and $H$.
Our main motivation in this section is to look for the possible functionals to remove the conditions (\ref{const}) on the radii $a$ and $r$  that were obtained by functionals depending only on the mean curvature. We shall show that there are some examples of such functionals where the torus is
a critical point without any constraint on radii $a$ and $r$.
Consider now the following functionals ${\cal F}_{n}$ defined by

\begin{equation}
{\cal F}^{K}_{n}=\int_{S}\, {\cal E}^{K}_{n}(H,K) dA+p \int_{V} dV
\end{equation}

\noindent where ${\cal E}^{K}_{n}$ are some functions of $H$ and $K$ for all $(n=1,2,3,\cdots)$, and $p$ is
a constant  and $V$ is the volume enclosed within the surface $S$. Below $a_{i},~ i=1-7$ are constants. We shall give two examples here and one more example in Appendix B.

\vspace{0.5cm}
\noindent
1. {\bf Third Order}:  ${\cal E}^{K}_{3}=a_{1}\,H^3+a_{2}\,H^2+a_{3}\,H+a_{4}+a_{5}\,K^2+a_{6}\,H\,K$

\vspace{0.5cm}
\noindent
{\bf Solution}:
\begin{eqnarray}
p&=&\frac{21 a_{1}-2a_{2}r-2a_{3}r^2}{2r^4},~~a_{6}=\frac{15a_{1}-2a_{2}r}{2}, \\
a_{5}&=&\frac{r(-15a_{1}+2a_{2}r)}{8},~~a_{4}=\frac{61a_{1}-6a_{2}r-8a_{3}r^2}{8r^3}
\end{eqnarray}
with $a^2=(6/5) r^2$. This is an example where the constraint $a^2 =(6/5) r^2$  at the third order is preserved with the inclusion
of $K$ terms into the Lagrange function ${\cal E}_{3}$. Adding also $H^2 K$ term to ${\cal E}_{3}$ does not change this fact.

\vspace{0.5cm}
\noindent
2. {\bf Fourth Order}: The most general fourth order functional is: \\
${\cal E}^{K}_{4}=a_{1}\,H^4+a_{2}\,H^3+a_{3}\,H^2+a_{4}\,H+a_{5}+a_{6}\,K^2+a_{7}\,H\,K+a_{8}\,H^2\,K$

\vspace{0.5cm}
\noindent
{\bf Solution}:
\begin{eqnarray}
p&=&\frac{1}{8r^5\,\Delta_{4}}[(-17a^6+71a^4 r^2+130a^2 r^4-120r^6)a_{1} -8\Delta_{4} r^2 (a_{3}+ r a_{4})],\\
a_{8}&=&\frac{a_{1}}{4 \Delta_{4}}\,(-55 a^6+23 a^4 r^2 -324 a^2 r^4 +144 r^6), \\
a_{7}&=&\frac{1}{8r \Delta_{4}}\,[(5 a^6+53 a^4 r^2+54 a^2 r^4-72 r^6)a_{1}-8 \Delta_{4} r^2 a_{3}],\\
a_{6}&=&\frac{1}{32 \Delta_{4}}\,[(270 r^6-1233 a^4 r^2 +1546 r^4 a^2 -648 r^6)a_{1}+8 r^2\,\Delta_{4} a_{3}], \\
a_{5}&=&\frac{1}{32 r^4 \Delta_{4}}\, [(-42 a^6+193 a^4 r^2 +39 a^2 r^4 -36 r^6) a_{1}-24 r^ 2 \Delta_{4} a_{3}-32 r^3 \Delta_{4} a_{4}],\\
a_{2}&=& \frac{a_{1}}{4 r \Delta_{4}}\, [-2 a^6-9 a^4 r^2 +38 a^2 r^4-24 r^6]
\end{eqnarray}
where
\begin{equation}
\Delta_{4}=(a^2-2 r^2)( a^2-r^2) (5 a^2-6 r^2)
\end{equation}
This example proves that as long as $\Delta_{5} \ne 0$ the radii $a$ and $r$ of the 2-torus are left arbitrary for functional ${\cal F}^{K}_{4}$ with
${\cal E}^{K}_{4}$ given above. If $\Delta_{4}=0$, since $a>r$, we have two choices $a^2=2 r^2$ and $a^2=(6/5)r^2$. Both of these choices lead to $a_{1}=0$ which reduces the degree of the polynomial ${\cal E}_{4}$ to 3 which is not acceptable.

\vspace{0.5cm}

\noindent
A special case where $a_{8}=0$:\\

${\cal E}^{K}_{4}=a_{1}\,H^4+a_{2}\,H^3+a_{3}\,H^2+a_{4}\,H+a_{5}+a_{6}\,K^2+a_{7}\,H\,K$

\vspace{0.5cm}
\noindent
{\bf Solution}:
\begin{eqnarray}
p&=&\frac{1168 a_{1}-5a_{3}r^2-5a_{4}r^3}{5r^5},~~a_{7}=\frac{783a_{1}-5a_{3}r^2}{5r}, ~~a_{2}=\frac{23a_{1}}{r},\\
a_{6}&=&\frac{r(-783a_{1}+5a_{3}r^2)}{20},~~a_{5}=\frac{3369a_{1}-15a_{3}r^2-20a_{4}r^3}{20r^4}
\end{eqnarray}
with $a^2=(12/11) r^2$. This is another example where the constraint (\ref{const4}) is preserved with the inclusion of some nontrivial $K$ terms.
We give one more example  ($n=5$) in Appendix B.

\section{ Main Theorems}

We start with the first theorem on the constraint equation in general.

\vspace{0.3cm}
\noindent
{\bf Theorem 1}:\,\,{\it 2-torus $T$ in ${\mathbb R}^3$ is a critical point of functionals ${\cal F}_{n}=\int_{S}\, {\cal E}_{n}\, dA+ p\, \int_{V}\,dV$ where ${\cal E}_{n}$ are the $n$th degree polynomials of the mean curvature $H$ of the surface $S$ if and only if
$\frac{a^2}{r^2}=\frac{n^2-n}{n^2-n-1}$ for all $n \ge 2$.}

\vspace{0.3cm}
\noindent
{\bf Proof}: Using Eqs.(\ref{tor1})-(\ref{tor3}) for the 2-torus we obtain

\begin{eqnarray}
\nabla^2 H&=&\frac{1}{a^2 r^3}\, [4r^3 (-a^2+r^2) H^3+2 r^2 (5 a^2 -6 r^2) H^2+4r(-2 a^2+3 r^2)H \nonumber \\
&&+2(a^2-2 r^2)], \\
(\nabla H) \cdot (\nabla H)&= &\frac{1}{a^2 r^4}\,[4 r^4 (-a^2+r^2)\, H^4+4 r^3 (3 a^2-4 r^2)\, H^3+r^2 (-13 a^2+24 r^2)\, H^2 \nonumber\\
&&+2 r (3 a^2-8 r^2)\, H+4 r^2]
\end{eqnarray}
Using these equations for the 2-torus  for all $n \ge 2 $ we get
\begin{equation}\label{identity1}
\nabla^2 H^n=\frac{4 n^2}{a^2}(-a^2 +r^2) H^{n+2}+\frac{2}{a^2 r}[(6 n^2-n)a^2-(8n^2-2n)r^2] H^{n+1}+\cdots
\end{equation}
Hence inserting
\begin{equation}
{\cal E}_{n}=a_{1}\,H^{n}+a_{2}\,H^{n-1}+\cdots+a_{n}\,H+a_{n+1},
\end{equation}\label{son1}
into the general shape equation (\ref{el1}) we get
\begin{eqnarray}
(n a_{1} \nabla^2 H^{n-1}+(n-1) a_{2} \nabla^2 H^{n-2}+\cdots)+(4 (n-1) a_{1} H^{n+1}+4(n-2)a_{2} H^{n}+ \cdots) \nonumber\\
-\frac{2}{r^2}(2 r H-1)( n a_{1} H^{n-1}+(n-1) a_{2} H^{n-2}+\cdots)=0. \nonumber
\end{eqnarray}
Using the identity (\ref{identity1}) for the 2-torus and collecting the coefficients of the powers of $H$ we get equations for $a_{i}$'s. The coefficient of the highest power $H^{n+1}$ in
(\ref{son1}) can be calculated exactly
\begin{equation}\label{sondenk}
a_{1}\,[\frac{4 n(n-1)^2\,(-a^2+r^2)}{a^2}+4(n-1)]\,H^{n+1}+\cdots=0.
\end{equation}
Then coefficient of $H^{n+1}$ must vanish which leads to the constraint equations

\[
\frac{a^2}{r^2}=\frac{n^2-n}{n^2-n-1},~~
\]
 for all $n \ge 2$. The remaining $n+1$ number of equations are linear algebraic equations for $a_{i}, ~~(i=1,2, \cdots, n+1)$ and $p$.
 In general one can solve them in terms of one arbitrary parameter, for instance $a_{1}$. In examples 1-6 in section 3 and Appendix A, the solutions contain two arbitrary coefficients. This means that one of the remaining equations in (\ref{sondenk}) vanishes identically. This completes the proof of Theorem 1.

 Our second theorem is on removing the constraint equations on the radii of the 2-torus by including $K$ terms into the Lagrange functions ${\cal E}_{n}$.

 \vspace{0.3cm}
\noindent
{\bf Theorem 2}:\,\,{\it 2-torus $T$ in ${\mathbb R}^3$ is a critical point of functionals ${\cal F}^{K}_{n}=\int_{S}\, {\cal E}^{K}_{n}\, dA+ p\, \int_{V}\,dV$ where ${\cal E}^{K}_{n}$ are appropriate polynomials of the mean and Gaussian curvatures  $(H,K)$ without any restriction
on the radii $a$ and $r$ }

\vspace{0.3cm}
\noindent
{\bf Proof}: We modify the Lagrange function ${\cal E}_{n}$ as
\begin{equation}
{\cal E}^{K}_{n}={\cal E}_{n}+\sum_{m=1}^{[n/2]}\, K^m \, \sum_{k=1}^{n-2m}\,a_{n+k+1} H^{k}
\end{equation}
where $[n/2]$ is $\frac{n}{2}$ for even $n$ and $\frac{n-1}{2}$ for odd $n$. There are additional $\frac{n^2}{4}$ parameters $a_{i}$ for even $n$ and $\frac{n^2-1}{4}$ for add $n$. The general shape equation (\ref{el1}) will pick up terms from the $K$
dependence hence we also need the following identities for 2-torus
\begin{eqnarray}
\nabla \cdot {\bar \nabla}\, H&=& \frac{1}{a^2 r^4}\,[16r^4(-a^2+r^2)\,H^4+12 r^3(4 a^2-5r^2)H^3 \nonumber \\
&&+4 r^2(-13 a^2+21 r^2)\,H^2+4 r (6 a^2-13 r^2)\,H \nonumber \\
&&+4(-a^2+3 r^2)],\\
\nabla \cdot {\bar \nabla}\, K&=& \frac{1}{a^2 r^5}\,[32 r^4 (-a^2+r^2)\,H^4+24 r^3 (4 a^2-5 r^2)\,H^3 \nonumber \\
&&+8 r^2 (-13 a^2+21 r^2)\,H^2+ 8 r (6 a^2-13 r^2)\,H \nonumber \\
&&+8(-a^2+3 r^2)]
\end{eqnarray}
There will be $n+1+\frac{n^2}{4}$ parameters $a_{i}$ in total for even $n$ and  $n+1+\frac{n^2-1}{4}$ for odd $n$. On the other the shape equation (\ref{el1}) leads to $n+1$ number of equations for coefficients $a_{i}$ and hence the coefficients of $H^{k}$ terms for $(k=0,1,2, \cdots, n, n+1)$ in (\ref{sondenk}) are modified and hence the constraint equation disappears. Since number equations are much less than the number of unknowns, modified Eq.(\ref{sondenk}) will be solved for $a_{i}$ leaving $\frac{n^2}{4}$ for even $n$ and $\frac{n^2-1}{4}$ for odd $n$  arbitrary. Since the coefficient of $H^{n+1}$ in (\ref{sondenk}) is modified with the new coefficients then constraint equation disappears.  Hence there remains a system of $n+1$ coupled linear equations for $n+1$ of coefficients $a_{i}$. Provided that determinant $\Delta_{n}$ of the coefficient matrix is different than zero, then  solutions of $a_{i}$  exists in terms of the arbitrary coefficients. If $\Delta_{n}=0$ the solution exists but with the constraint equation. This completes the proof of Theorem2. In section 4 and Appendix B we give two examples $n=4$ and $n=5$ for illustration of this theorem.

\section{Conclusion}

Clifford torus is a critical point of  Willmore and also Helfrich's functional where the Lagrange function is a quadratic polynomial
in the mean curvature of a closed surface in ${\mathbb R}^3$. The main contribution of this work is that the sequence of tori
$\{{T}_{n}\}$ where radii $a$ and $r$ restricted to satisfy  $\frac{a^2}{r^2}=\frac{n^2-n}{n^2-n-1}$ for all $n \ge 2$ are the critical
points of the functionals ${\cal F}_{n}$ where the Lagrange function ${\cal E}_{n}$ are polynomial of degree $n$  in the mean curvature $H$ of the surface. We have given 6 examples $n=1-6$ and proved this assertion in section 5. If the Lagrange functions ${\cal E}_{n}$ are modified by adding
some $K$ terms, Gaussian curvature of the surfaces, it is possible to remove the constraint equations. We have given examples for $n=4$ and $n=5$.
We have proved in general that 2-torus is a critical point of sequence of functionals ${\cal F}_{n}$ where the Lagrange functions ${\cal E}_{n}$
is a polynomial of both $H$ and $K$.

In section 3 and Appendix A, for each solution with $n=2-6$ and $p=0$ we have calculated the curvature energy ${\cal F}_{n}$. As in the case of the Willmore energy functional ($n=2$) it is expected that the tori for $n \ge 3$  with the constraints, are minimal energy surfaces. To support this assertion, second variation of the functionals on these surfaces must be calculated. Another point to be examined is the stability of these minimal energy surfaces. These points will be communicated  in a forthcoming publication.

\section{Appendix A}

In this Appendix we add three more examples to those given in section 3.

\vspace{0.3cm}
\noindent
4. {\bf Fourth Order Functional}:  ${\cal E}_{4}=a_{1}\,H^4+a_{2}\,H^3+a_{3}\,H^2+a_{4}\,H+a_{5}$.

\vspace{0.5cm}
\noindent
{\bf Solution}:
\begin{eqnarray}
p&=&\frac{77 a_{1}- a_{4} r^3}{ r^3},~~a_{5}=\frac{51 a_{1}- a_{4} r^3}{r^4},~~a_{3}=\frac{783 a_{1}}{5 r^2},~~a_{2}=\frac{23 a_{1}}{r} ,\\ a^2&=&(12/11) r^2 \label{const4}
\end{eqnarray}
We find that (for $p=0$)
\begin{equation}\label{min-en4}
{\cal F}_{4}=\int_{T}\, {\cal E}_{4} dA=\frac{666\, \sqrt{11} \,\pi^2}{5} ( \frac{a{1}}{r^2})
\end{equation}

\vspace{0.5cm}
\noindent
5. {\bf Fifth Order Functional}:  ${\cal E}_{5}=a_{1}\,H^5+a_{2}\,H^4+a_{3}\,H^3+a_{4}\,H^2+a_{5}\,H+a_{6}$.

\vspace{0.5cm}
\noindent
{\bf Solution}:
\begin{eqnarray}
p&=&\frac{41915\, a_{1}-14\, a_{5} r^4}{14 r^6},~~~ a_{6}=\frac{13848 a_{1}-7 a_{5} r^4}{7 r^5},\\
a_{4}&=&\frac{715045 a_{1}}{126 r^3}, ~~~a_{3}=\frac{12035 a_{1}}{14 r^2},~~~a_{2}=\frac{50 a_{1}}{r},\\
a^2&=&(20/19) r^2
\end{eqnarray}
we find that (for $p=0$)

\begin{equation}
{\cal F}_{5}=\int_{T}\, {\cal E}_{5} dA=\frac{235750 \sqrt{19}}{63}\, \pi^2\, (a_{1}/r^3)
\end{equation}

\vspace{0.5cm}
\noindent
6. {\bf Sixth Order Functional}:  ${\cal E}_{6}=a_{1}\,H^6+a_{2}\,H^5+a_{3}\,H^4+a_{4}\,H^3+a_{5}\,H^2+a_{6}\,H+a_{7}$.

\vspace{0.5cm}
\noindent
{\bf Solution}:
\begin{eqnarray}
p&=&\frac{544813 a_{1}-32 a_{6} r^5}{32 r^7},~~a_{7}=\frac{1796815 a_{1}-16 a_{6} r^5}{16 r^6},\\
a_{5}&=&\frac{139780065 a_{1}}{448 r^4},~~~a_{4}=\frac{1534015 a_{1}}{32 r^3}, \\
a_{3}&=&\frac{6235 a_{1}}{2 r^2},~~~a_{2}=\frac{183 a_{1}}{2 r}, \\
a^2&=&(30/29) r^2
\end{eqnarray}
we find that (for $p=0$)

\begin{equation}
{\cal F}_{4}=\int_{T}\, {\cal E}_{6} dA=\frac{37643625 \sqrt{29}}{224}\, \pi^2\, (a_{1}/r^4)
\end{equation}

\section{Appendix B}

In this Appendix we add one more example those given in section 4.

\vspace{0.3cm}
\noindent
3. {\bf Fifth Order}: The most general fifth order functional is: \\
${\cal E}^{K}_{4}=a_{1}\,H^5+a_{2}\,H^4+a_{3}\,H^3+a_{4}\,H^2+a_{5}\,H+a_{6}+a_{7}\,K^2\,H+a_{8}\,K^2+a_{9}\,K\,H^3\,+a_{10} K H^2+a_{11}\,K H$

\vspace{0.5cm}
\noindent
{\bf Solution}:
\begin{eqnarray}
p&=&\frac{1}{96 r^6\,\Delta_{5}      }[(-150a^8+1174a^6 r^2+7218a^2 r^6-2520r^8)a_{1} \nonumber \\
&&+r\,(-240 a^8+930 a^4 r^2+1554 a^4 r^4 -4116 a^2 r^7+1872 r^8)\,a_{2}-72 r^2 \Delta_{5}\, a_{3} \nonumber \\
&&-480 r^3 \Delta_{5}\,a_{4}-96 r^4 \Delta_{5}\, a_{5}],\\
a_{11}&=&\frac{1}{48 r^2 \Delta_{5}}\,[(66 a^8-721 a^6 r^2-4 a^4 r^4+924 a^2 r^6 -480 r^8)\,a_{1} \nonumber \\
&&+r\,(18 a^8+246 a^6 r^2+288 a^4 r^4+1128 a^2 r^6+576 r^8)\,a_{2} \nonumber \\
&&-24 r^2 \Delta_{5} (a_{3}+2r a_{4})], \\
a_{10}&=&\frac{1}{24 r}\,[(226 a^8-1141 a^6 r^2+2161 a^4 r^4 -1866 a^2 r^6+600 r^8)\, a_{1} \nonumber \\
&&+r\,(-342 a^8 +1704 a^6 r^2-3700 a^4 r^4 +2436 a^2 r^4-720 r^8)\, a_{2} \nonumber \\
&&-24 r^2 \Delta_{5}\, a_{3}], \\
a_{9}&=&\frac{1}{6 \Delta_{5}}\,(a^2-r^2)(a^2-2 r^2)(a^2-\frac{20}{19} r^2)(95 a^2-114 r^2),\\
a_{8}&=&\frac{1}{192 r \Delta_{5}}\,(-2702 a^8+11838 a^6 r^2-25525 a^4 r^4 +18894 a^2 r^6-5880 r^8)a_{1} \nonumber \\
&&+r\,(1644 r^8-8934 a^6 r^2+16230 a^4 r^4-12540 a^2 r^6 +3600 r^8)\, a_{2} \nonumber \\
&&+48 r^3 \Delta_{5} (a_{3}+r a_{4})],\\
a_{7}&=&\frac{1}{32 \Delta_{5}}\,[(358 a^8-1948 a^6 r^2+3833 a^4 r^4-3238 a^2 r^+1000 r^8)a_{1} \nonumber\\
&&+r\,(4 a^8+14 a^6 r^2 -94 a^4 r^4 +124 a^2 r^6 -48 r^8)a_{2}+8 \Delta_{5} r^2 a_{3}],\\
a_{6}&=&\frac{1}{64 r^5 \Delta_{5}}\,[(-60 a^8+458 a^6 r^2 -2913 r^4 a^4 +3398 a^2 r^6-1208 r^8)a_{1}\nonumber \\
&&+r\,(-100 a^8+418 a^6 r^2 +766 a^4 r^4-1996 a^2 r^6+912 r^8)a_{2} \nonumber \\
&&-32 r^2 \Delta_{5} a_{3}-48 r^3 \Delta_{5} a_{4}-64 r^4 \Delta_{5} a_{5}],
\end{eqnarray}
where
\begin{equation}
\Delta_{5}=(a^2- r^2)^2\,(a^2-2 r^2) (5 a^2-6 r^2).
\end{equation}
Again as long as $\Delta_{5} \ne 0$ we have here another example of a functional where an arbitrary 2-torus is a critical point. If $\Delta_{5}=0$,
since $a>r$, we have only two cases:

\vspace{0.3cm}
\noindent
{\bf Case 1}: $a^2=2 r^2$.

\begin{eqnarray}
p&=&\frac{1}{108 r^6}\,[-31 a_{1}-9 a_{3} r^2+36 a_{5} r^4+144 a_{6} r^5], \\
a_{11}&=& \frac{1}{54 r^2}\,[97 a_{1}+9 a_{3} r^2+72 a_{5} r^4+72 a_{6} r^5], \\
a_{10}&=&-\frac{1}{9r}\,[38 a_{1}+9 a_{3} r^2],~~~a_{9}=-3a_{1},\\
a_{8}&=&\frac{1}{108 r}\,[163 a_{1}+9 a_{3} r^2-36 a_{5} r^4-36 a_{6} r^5], \\
a_{7}&=&\frac{1}{36}\,[73 a_{1}+9 a_{3} r^2],~~~a_{2}=\frac{35 a_{1}}{18 r}.
\end{eqnarray}
The coefficients $a_{1}$, $a_{3}$, $a_{4}$, $a_{5}$ and  $a_{6}$ are left arbitrary.

\vspace{0.3cm}
\noindent
{\bf Case 2}: $a^2=\frac{6}{5} r^2$.

\begin{eqnarray}
p&=&\frac{1}{8424 r^6}\,[-9119 a_{1}-486 a_{3} r^2+324 a_{4} r^3+3240 a_{5}r^4+11664 a_{6} r^5], \\
a_{11}&=& \frac{1}{21060 r^2}\,[88723 a_{1}-162 a_{3} r^2+-5508 a_{4} r^3+20736 a_{5} r^4+20736 a_{6} r^5], \\
a_{10}&=&-\frac{1}{10530r}\,[-155972 a_{1}-9882 a_{3} r^2+972 a_{4} R^3+296 a_{5} r^4 \nonumber \\
&&+1296 a_{6} r^5],\\
a_{8}&=&\frac{1}{42120 r}\,[377387 a_{1}+5022 a_{3} r^2+2268 a_{4} r^3-11016 a_{5} r^4 \nonumber \\
&&-11016 a_{6} r^5], \\
a_{7}&=&\frac{1}{42120}\,[39317 a_{1}+9882 a_{3} r^2-912 a_{4} r^3-1296 a_{5} r^4-1296 a_{6} r^5],\\
a_{9}&=&-\frac{7a_{1}}{5},~~~a_{2}=\frac{791 a_{1}}{162 r}.
\end{eqnarray}
The coefficients $a_{1}$, $a_{3}$, $a_{5}$ and  $a_{6}$ are left arbitrary.

We claim that
for all functionals with ${\cal E}_{n}=a_{1}\, H^n+\cdots $ discussed in section 3, it is possible to add suitable $K$ terms so that an arbitrary 2-torus is a critical point
without any constraint on the radii $a$ and $r$. We prove our claim in  section 5.

\section{Acknowledgments}
 This work is partially supported by the Scientific and Technological
Research Council of Turkey (TUBITAK).

\end{document}